\documentclass[11pt]{article}
\usepackage{mathptmx}
\usepackage{latexsym,amssymb,amsmath,amscd,amscd,amsthm,amsxtra}
\usepackage[dvips]{graphicx}

\newtheorem{thm}{Theorem}[section]
\newtheorem{lem}[thm]{Lemma}

\newtheorem{pro}[thm]{Proposition}

\newcommand{\px}{\frac{\partial }{\partial x}}
\newcommand{\py}{\frac{\partial }{\partial y}}
\newcommand{\pz}{\frac{\partial }{\partial z}}

\newcommand{\lon }{\,\rightarrow\,}

\newcommand{\pf}{\noindent{\bf Proof.}\ }

\newcommand{\frkg}{\mathfrak g}

\def\gpd{\,\lower1pt\hbox{$\longrightarrow$}\hskip-.24in\raise2pt
         \hbox{$\longrightarrow$}\,}

\def\qed{\hfill ~\vrule height6pt width6pt depth0pt}

\begin{document}
\title{{Quadratic Deformations of Lie-Poisson Structures
\thanks{ Research partially supported by NSF of China  and the Research Project of ``Nonlinear Science".}}}
\author{ Qian LIN, ~ Zhangju LIU ~and ~Yunhe SHENG \\
   Department of Mathematics and LMAM\\ Peking University, ~~
Beijing 100871, China\\
          {\sf email: linqian@pku.edu.cn;\quad liuzj@pku.edu.cn;\quad syh@math.pku.edu.cn} }
\date{}
\maketitle
\begin{abstract}
In this letter, first we give a decomposition  for any Lie-Poisson
structure $\pi_{\frkg}$   associated to the modular vector. In
particular, $\pi_{\frkg}$ splits into
         two compatible   Lie-Poisson structures if  $dim{\frkg} \leq 3$.
        As an application, we classified quadratic deformations of Lie-Poisson structures on $\mathbb
        R^3$ up to linear diffeomorphisms.
\end{abstract}
\section{Introduction}

It is known that linear  and quadratic Poisson structures are two
most basic and important Poisson structures both for their rich
algebraic and geometric properties and  various applications in
physics and other fields of mathematics.
 The linear Poisson structures are in one-to-one correspondence
with  Lie algebra structures and usually called Lie-Poisson
structures.  The idea of using linear Poisson brackets to
understand the structure of Lie algebras can be traced back to the
work of Lie. In this spirit there have been some suggestions of
pursuing this geometric approach for Lie algebra structures( e.g.,
see \cite{CIMP}, \cite{GMP} and \cite{Sheng}).  Besides the linear
Poisson structures,   quadratic Poisson structures are also
generally studied from various aspects (e.g., \cite{DM},
\cite{DH}, \cite{MMR} and see \cite{We2 } for more comments).

 The  purpose of  this letter is to study quadratic deformations of
      Lie-Poisson structures  and  their classification on $\mathbb{R}^3$.
      Main motivation for us comes from the work in
      \cite{CS} and \cite{CMS }, where some special quadratic deformation of a
      Lie-Poisson structure appears when the authors study some
      geometric objects such as holonomy and symplectic
      connection.  In \cite{CF} a simple example is also shown as  an example of Poisson-Dirac submanifolds.
       On the other hand, in \cite{DH} and \cite{LX},       the
     quadratic Poisson structures on $\mathbb{R}^3$, which    can be considered as quadratic
deformations of the abelian Lie-Poisson structure, are totally
classified. Thus, it is natural  to consider  all possible
quadratic deformations for any Lie-Poisson structure. To classify
quadratic deformations for a fixed Lie-Poisson structure it is
enough  to use linear transformations, which keep the degree of  a
homogenous tensor field.
  The classification under general diffeomorphism is more
  complicated. Please see \cite{DZ} for the classification under local diffeomorphisms of every 3-dim
Poisson structure vanishing at a point with a non-zero linear
part.

To save space we just write out details of quadratic deformations
for some spacial cases. The others can be done by same way without
any difficulty except for some tedious computations.

{\bf Acknowledgement:} ~~ The second author would like to thank
Professor T. Ratiu for  useful comments and hospitality during his
visit at the Bernoulli Center. We are also grateful to G. Marmo
and the referee  for many  useful comments and suggestions.
\section{The Classification of Lie-Poisson structures on $\mathbb R^3$}
Let  $\Omega = dx_{1}\wedge dx_{2}\cdots\wedge dx_{n}$ be the
canonical volume
        form  on $\mathbb{R}^n$. Then  $\Omega$ induces an
isomorphism  $\Phi $ from the space of all $i$-multiple vector
fields to the space of all $(n-i)$-forms. Let $d$ denote the usual
exterior differential on forms and
 $$D=(-1)^{k+1}\Phi^{-1} \circ d \circ  \Phi: ~~\mathcal{X}^{k}(\mathbb{R}^n) \lon
\mathcal{X}^{k-1}(\mathbb{R}^n),$$ its pull back under the
isomorphism $\Phi$. The Schouten bracket can be written in terms
of this operator as follows\cite{Ko}:
\begin{equation}
\label{eq:schouten} [U,\ V]=D(U\wedge V)-D(U)\wedge V-(-1)^{i}
U\wedge D(V),
\end{equation}
for all $U\in \mathcal{X}^{i}(\mathbb{R}^n)$ and $V \in
\mathcal{X}^{j}(\mathbb{R}^n)$.
       It is obvious that  there is a one-to-one
correspondence between matrices in $\mathfrak{gl(n)}$ and linear
vector fields on $\mathbb{R}^n$, i.e.,
\begin{equation}\label{eq Ahat} A =(a_{ij}) \longleftrightarrow \hat{A}=\sum_{ij}a_{ij}x_{j}\frac{\partial}{\partial
x_{i}} , ~~~~~~ \, ~~~ div_{\Omega}\hat{A}= D(\hat{A}) = \mbox{tr}
A.
\end{equation}
 A vector $k \in \mathbb{R}^n$ corresponds a
constant vector field $\hat{k}$ by translation and satisfies
\begin{equation}\label{eq Ak}
div_{\Omega}\hat{k}= D(\hat{k}) = 0, \, ~~~~~~~~~\,
~~~~~~[\hat{A},~\hat{k} ] = -\hat{Ak},~~~~~~~~~\,~~~~~~\forall ~A
\in \mathfrak{gl(n)}.
\end{equation}

 For  a given  Poisson tensor $\pi$, let
$D(\pi )$ be its modular vector field (see \cite{We}, which is also
called the curl vector field in \cite{DH}). Such a vector field
 is always compatible with $\pi$, i.e.,
$[D(\pi),\ \pi]=0$.
 A Poisson  structure is called unimodular if
$ D(\pi ) = 0$.  By (\ref{eq Ak}), it is easy to see that, for any
$k \in \mathbb{R}^n$,  The bi-vector field $\hat{I}\wedge\hat{k}$
is a linear Poisson structure with the modular vector $(n-1)k$,
where $I$ is the identity matrix.
 The corresponding Lie algebra is called {\em book
 algebra} when
 $ n =3$ and  $k = (0, 0, 1)$.
 Next we give a  similar
        decomposition for  Lie-Poisson structures as doing in  \cite{LX}  for quadratic Poisson structures.
\begin{thm}\label{Th. decom}
            Any linear Poisson structure $\pi$ on $  \mathbb{R}^n $ has a unique
            decomposition:\\
            \begin{equation}\label{linear}
            \pi = \frac{1}{n-1}\hat{I}\wedge\hat{k} + \Lambda,
            \end{equation}
            where $k \in  \mathbb R^n$ is
            the modular vector of $\pi$ and $\Lambda$ is a linear bi-vector
            field such that
 $D(\Lambda)=0$. $\Lambda$ is a Poisson   structure  compatible with $\pi$ if and only if $D(\Lambda
\wedge \Lambda)=0. $ In particular, $\Lambda$ is always a
unimodular Lie-Poisson structure  on  $ \mathbb R^3$.
        \end{thm}
        \pf Let $\hat{k} = D(\pi)$ and
        define a bi-vector field
        $\Lambda =:\pi-\frac{1}{n-1}\hat{I}\wedge\hat{k}$.
        Then
       $D(\Lambda)=0$  and   $ [ \Lambda,\Lambda ] = D(\Lambda
\wedge \Lambda)$ by   formula (\ref{eq:schouten}).  The fact that
$[\hat{k} ,\Lambda ] = 0$  is because $\hat{k}$ is the  modular
vector field for both Poisson structures $\pi$ and $
\hat{I}\wedge\hat{k}$. Finally, by  (\ref{eq:schouten}) and the
fact that $[\hat{I}, \pi] = -\pi$, one can see that $\Lambda$ is
compatible with $\pi$ if  $\Lambda$  is a Poisson   structure.
\qed

From now on, we  focus on the case that $n=3$. In this case
  there
exists  a  unique homogeneous quadratic   3-vector field $L$ such
that $\Lambda=D(L)$ since $ D(\Lambda) = 0$ and the cohomology
groups here are trivial. The quadratic 3-vector fields are in
one-to-one correspondence with quadratic functions $f$ via the
volume form $\Omega$. We use $\pi_{f}$ to denote this unimodular
linear
 Poisson structure  which is given  by
\begin{equation}\label{eq exact}
\pi_{f}=\Phi^{-1}(df)= \frac{\partial f}{\partial x}\py \wedge \pz+
\frac{\partial f}{\partial y}\pz \wedge \px+ \frac{\partial
f}{\partial z}\px \wedge \py.
\end{equation}
The compatibility condition with the modular vector field $\hat{k}$
is given by
\begin{equation}\label{eq kf}
[ \hat{k} , \pi_{f}] = \Phi^{-1}L_{\hat{k}}(df) =
\Phi^{-1}d(\hat{k}f) =0 \, ~~~~~~~~ \Longleftrightarrow \,
~~~~~~~~~~~ \hat{k}f = 0.
\end{equation}
 Therefore, the linear Poisson structures   on    $\mathbb{R}^3$   are in one-to-one correspondence with  pairs $(k,f)$,
 where $k$ is a
 vector and $f$ is a quadratic function, such that $\hat{k}f = 0$.
As doing in \cite{LX} for the quadratic Poisson structures, here we
also call such a pair $(k, f)$ as a {\bf{compatible
        pair}} of the corresponding
        linear Poisson structure $\pi$ and usually denote $\pi=\pi_{k,f}$. The next result characterizes the isomorphisms of linear Poisson structures by
means of their compatible pairs.
\begin{thm}\label{Th. Isomor}
            Let $\pi_{1}$ and $\pi_{2}$ be two linear Poisson structures on $\mathbb{R}^3$ determined
            by  compatible pairs
            $( k_{1},f_{1} )$ and $( k_{2},f_{2} )$ respectively. Then $\pi_{1}$ is
            isomorphic to $\pi_{2}$ if and only if there is a $T \in GL(3)$ such that
                         $$k_{2}=Tk_{1},\qquad f_{2}=det(T)f_{1}\circ T^{-1}.$$
Particularly, the automorphism group of $\pi_{k,f}$, denoted by
$Aut(\pi_{k,f})$, is
\begin{equation}
                Aut(\pi_{k,f}) = \{T |~~ T \in GL(3),~~ Tk = k,~~f\circ T =
                det(T)f\}.
\end{equation}
Consequently, for the corresponding Lie algebra $\frak g$ of
$\pi_{k,f}$, one has
\begin{equation}
              Der(\frkg) \cong \{D|~~ D \in {\frak gl(3)},~~ Dk = 0 ,~~  \hat{D}f
              =(trD)f\}.
\end{equation}

 \end{thm}
        \pf $T$ is an isomorphism means that $T_{*}\pi_{1}=\pi_{2}$, by
        Theorem \ref{Th. decom} and  properties of the modular
        vector, which is equivalent to that $Tk_{1}=k_{2}$ and
        $T_{*}\pi_{f_{1}}=\pi_{f_{2}}$, which is equivalent to
       $|T|f_{1} = f_{2}\circ T$. The other conclusions are
easy to be checked.\qed

 Moreover, it will be seen that the traceless derivations,
 \begin{equation}\label{eq tr0}
               Der_0(\frkg) \cong \{D |~~ D \in {\frak sl(3)},~~ Dk = 0 ,~~  \hat{D}f = 0\},
 \end{equation}
play an important role for the quadratic deformation. As an
application of Theorems \ref{Th. decom} and
 \ref{Th. Isomor}, we show  a simple way to classify linear Poisson structures  on
 $\mathbb{R}^3$.

\begin{thm}\label{Th classify}
            Any linear Poisson structure $\pi_{k, f}$ on $\mathbb{R}^3$ is isomorphic to one of
            the following standard forms:

                  \setlength{\parindent}{15mm} {\rm(A)}.~ $k=0$ (unimodular )\qquad\qquad{\rm(B)}.~$k=(0,0,1)^{T}$,~~ \, $\hat{k}= \frac{\partial}{\partial z}$,
                  \begin{align*}
         (1).~f&=0         &     (7). ~ f&=0\\
         (2).~f&=x^2+y^2+z^2 &   (8).~f&=a(x^2+y^2), \, \, a> 0\\
          (3).~ f&=x^2+y^2-z^2 & (9).~f&=a(x^2-y^2), \, \, a> 0\\
           (4).~ f&=x^2+y^2 &    (10).~f&=x^2.\\
         (5).~  f&=x^2-y^2 &\\
             (6).~ f&=x^2&
\end{align*}
      \end{thm}
        \pf For  unimodular cases, $\pi_{f_1}$ is isomorphic to $\pi_{f_2}$ if and only if there is some $T \in GL(3)$ such that
$f_2\circ T=\mid T\mid f_1$ by Theorem \ref{Th. Isomor}.  First we
take a  $T^\prime \in GL(3)$ such that the quadratic function
$f^\prime = f\circ T^\prime$ is
 one of the standard   forms  listed from $(1)$ to $(6).$
     Assume $T=det(T^\prime)^{-1}T^\prime$, then we
        have $f\circ T=\mid T\mid f^\prime.$

               If $k\neq0$, we can take $k=(0,0,1)^{T}$ by a coordinate transformation.
               Then $\exists A \in symm(2)$ such that
               $f = (x, y)A(x, y)^T$
               because
        $\hat{k}f = \frac{\partial f}{\partial z}= 0$ by (\ref{eq kf}). Moreover, $T$ must has the
        form:
        $ \left(
            \begin{array}
            {cc}
            S         & 0\\
            \alpha    & 1
            \end{array}
        \right),$  where $ S \in GL(2),~\alpha\in \mathbb{R}^2$,  since $Tk = k$ by Theorem \ref{Th. Isomor}.
        Thus it is seen that,  for any $T$  having the form above and $f \neq 0$, the induced new quadratic function from
        $f$ must be determined by the matrix $|S|^{-1}S^TAS$ so that its standard form can be  fixed up  to a
        constant. Furthermore it is easy to check that this constant can be adjust by a sign and have the form that listed
        in the theorem. As in  case (10), we can adjust $S$ such that $det(S)=1$.\qed

Denote $\frak g_i(i=1,\cdots,10)$ the corresponding Lie algebra of
linear Poisson structures $\pi_i$.   Comparing with the
classification described in \cite{Jac} via the
     dimension  of  the derived algebra
       $[\frak g,\frak
    g]$,  Case (1) is abelian. Cases (2) and
    (3) are
 simple Lie algebras.
In cases (4), (5), (7)-(10) the derived algebra has dimension 2
except for $a=1/4$ in  case (9). In  case (6) and  case (9) with
$a=1/4$, the derived algebra has dimension 1.
    The Lie algebras (8)-(10) are three twisted Lie
    algebras
    of $\frak g_7$ corresponding
    three 2-cocycles respectively.

   By means of Theorem \ref{Th. Isomor},  we see that the automorphism
   groups are independent on the constant in Cases (8)-(9).
   Combining with
 Theorem \ref{Th classify},  we can get the automorphism
groups of  3-dim. linear Poisson structures.

 \begin{thm}\label{Th Auto}
        Let $G_i, ~ i = 1, \cdots, 10$ denote the automorphism group of the linear Poisson structure which corresponds to Case (i) in Theorem \ref{Th classify}.
         Then we have
               \begin{itemize}
            \item[] $G_1 = GL(3)$, $G_2 = SO(3)$, $G_3 = SO(2,1)$, $G_4 = \{\left(
                                        \begin{array}
                                        {cc}
                                        \lambda T   &   0      \\
                                        \xi      &  |T|
                                        \end{array}
                                \right)| ~~T \in
                                O(2)\}$,
            \item[] $G_5 = \{\left(
                                        \begin{array}
                                        {ccc}
                                        \lambda\alpha      &   \beta        &   0 \\
                                        \lambda\beta       &   \alpha       &   0 \\
                                        \gamma      &   \delta       &  \lambda

                                    \end{array}
                                \right) |~~ \lambda = \pm1,
                                 ~~\alpha^2\neq\beta^2\}$,
            ~~ \, ~~  $G_6 = \{\left(
                                        \begin{array}
                                        {cc}
                                        a           &   0     \\
                                        \xi       &   A

                                    \end{array}
                                \right)|
                                 ~~|A| = a \}$,
            \item[]  $G_7 = \{\left(
                                        \begin{array}
                                        {cc}
                                        A       &      0 \\
                                        \xi &      1

                                    \end{array}
                                \right)|~~A\in GL(2)\}$,
            ~~ \, ~~ $G_8 = \{\left(
                                        \begin{array}
                                        {cc}
                                        \lambda T   &   0      \\
                                        \xi       &   1

                                        \end{array}
                                \right)| ~~T \in SO(2), ~~\lambda \neq 0\}$,
            \item[]  $G_9 = \{\left(
                                        \begin{array}
                                        {ccc}
                                        \alpha      &   \beta        &   0 \\
                                        \beta       &   \alpha       &   0 \\
                                        \gamma      &   \delta       &   1

                                    \end{array}
                                \right) |~~ \alpha^2\neq\beta^2\}$,
            ~~ \, ~~  $G_{10} = \{\left(
                                        \begin{array}
                                        {ccc}
                                        \alpha      &   0            &   0 \\
                                        \beta       &   \alpha       &   0 \\
                                        \gamma      &   \delta       &   1

                                    \end{array}
                                \right) |~~\alpha\neq
                                0\}$.
        \end{itemize}
    \end{thm}

Note that $G_{8}$-$G_{10}$  are subgroups of $G_{4}$-$G_{6}$
respectively with $Tk = k$  and    subgroups of $G_{7}$. This fact
will be used  in the last section to classify quadratic
deformations of Lie-Poisson structures  of Cases (8)-(10).

\section{ Quadratic Deformations}

In \cite{LX}, any
     quadratic Poisson structure on $\mathbb{R}^3$ is characterized by its compatible
     pair $(K,F)$, where $K \in {\frak sl(3)}$ and $F$
     is a homogeneous cubic polynomial such that $\hat{K}F =0$ and
     $\pi_{_K,_F} = \pi_{_F}+\frac{1}{3}\hat{I}\wedge\hat{K}.
$ Such a Poisson structure can be considered as a quadratic
deformation of the abelian Lie-Poisson structure.

For any Lie-Poisson structure $\pi_{k,f}$, we shall  study its
quadratic deformations $\pi_{k,f}+\pi_{_K,_F}$,  evidently,  which
is still a Poisson structure if and only if
$[\pi_{k,f},\pi_{_K,_F}]=0$.
For convenience, denote by $\mathcal{C}$ the
  compatible pairs of the quadratic Poisson structures and  $\mathcal{C}_{k,f}$ the compatible pairs that can make
quadratic deformations of the linear Poisson structure
$\pi_{k,f}$, i.e.,
\begin{equation}\label{eq C}
\mathcal{C}=\{(K,F);~K\in {\frak sl(3)}, ~~  F
     ~\mbox{is a 3- polynomial ~s.t.} ~\hat{K}F =0\}.
\end{equation}
\begin{equation}\label{eq C kf}
\mathcal{C}_{k,f}=\{(K,F)\in\mathcal{C},~[\pi_{k,f},\pi_{_K,_F}]=0\}.
\end{equation}
Denote $X=(x,y,z) \in \mathbb R^3$ and for any quadratic function
$f$, write $f=(AX,X)$, where $A$ is a symmetric matrix. For any
$k\in\frak g^*$, denote $\widetilde{k}:\frak g^*\rightarrow \frak
g$  the
          skew-symmetric matrix corresponding to $\Phi^{-1}(k) \in
          \frak g \wedge \frak g$. Then we have

     \begin{thm}\label{Th deform} Let $(k,f)$ and $(K,F)$ be two  compatible
     pairs of a Lie-Poisson structure and a quadratic Poisson structure on $\mathbb{R}^3$  respectively.
     Then we have
        \begin{equation}\label{eq 1}
         (K,F)\in \mathcal C_{k,f}\Longleftrightarrow
        \widehat{k}F=-\frac{1}{6}X(12A+\widetilde{k})KX^T.
         \end{equation}
Especially, in   unimodular cases that $k=0$, we have
\begin{equation}
\mathcal C_{0,f}=\{(K,F)\in\mathcal C,~\widehat{K}f=0\},~~
i.e.,~~~~ K\in Der_0({\frak g}).
\end{equation}
 \end{thm}
     \pf By means of Equality (\ref{eq:schouten}), $\pi_{_K,_F}$  makes quadratic deformation
     of
     Lie-Poisson struture $\pi_{k,f}$ if and only if
$$-[\pi_{k,f},\pi_{_K,_F}] =
\widehat{k}\wedge(\frac{1}{3}\widehat{I}\wedge\widehat{K}+\pi_{_F})+\widehat{K}\wedge(\frac{1}{2}\widehat{I}\wedge\widehat{k}+\pi_f)=0.
$$
This is just 6$
\pi_{\widehat{k}F+\widehat{K}f}=D(\widehat{K}\wedge\widehat{k}\wedge\widehat{I})
$ and is equivalent to
$$
6(\widehat{k}F+\widehat{K}f)=-(x,y,z)(\widetilde{k}K)(x,y,z)^T.
$$
Note that
$\widehat{K}f=(x,y,z)(AK+(AK)^T)(x,y,z)^T$, 
so a compatible pair $(K,F)\in \mathcal C_{k,f}$  if and only if
it satisfies Equality (\ref{eq 1}). \qed

Next we consider the problem to classify the quadratic
deformations of a Lie-Poisson structure on $\mathbb{R}^3$. We
restrict us to the case that the linear parts of two isomorphic
Poisson structures are also isomorphic. That is, we consider the
classification of  quadratic deformations for a fixed Lie-Poisson
structure by using linear transformations, which keep the degree
of  a homogenous tensor field.  The following theorem is
straightforward.
            \begin{thm}\label{Th 2iso} Let $(K_i,F_i)_{i=1,2}\in \mathcal
            C_{k,f}$,
                where $(k,f)$ is the compatible pair of Lie-Poisson  structure $\pi_{k,f}$.
                Then $\pi_{k,f}+\pi_{K_i,F_i}$ are isomorphic if and only if there is  a $T \in Aut(\pi_{k,f})$
                such that
                \begin{equation}\label{lx}
                K_{2}=TK_{1}T^{-1} ~~ and ~~ F_{2}=det(T)F_{1}\circ T^{-1}.
                \end{equation}
                Consequently, the classification of quadratic deformations of a Lie-Poisson
structure $\pi_{k,f}$ on $\mathbb R^3$ is  parameterized  by the
orbit space $ Aut(\pi_{k,f}) \backslash \mathcal{C}_{k,f}$.
\end{thm}

 It is
known that $GL(3)$ acts on the space of  compatible pairs
$\mathcal{C}$ given in (\ref{eq C}) as
\begin{equation}
\label{eq:pair} T: \mathcal{C} \rightarrow \mathcal{C}, ~~~ T(K,F)
=: (TKT^{-1},\, det(T)F\circ T^{-1}),  ~~~~  ~~ \forall~T \in
GL(3),
\end{equation}
so that, as did in \cite{LX}, the quadratic Poisson structures on
$\mathbb{R}^3$ are classified by the Adjoint  orbits of $GL(3)$
and one can take the Jordan forms as the standard forms. The above
theorem shows that one should do more things on the orbits of the
Jordan forms to classify the quadratic deformations of a
Lie-Poisson structure.

Notice that $\mathcal{C}_{k,f}$ is relatively easy to be
determined in the unimodular cases but it is  difficult in Cases
$(7)-(10)$ because the unknown data $K$ and $F$ are involved
together in Equation (\ref{eq 1}).  In these cases we should fix
$K$ firstly with some standard form and then to find all possible
compatible 3-polynomials $F$  satisfying Equation (\ref{eq 1}).
For a fixed Jordan form $K\in\frak sl(3)$, denote by $G_K \subset
GL(3)$ as its isotropy subgroup for the adjoint action and
$\mathcal J_K= GL(3)/G_K$ the adjoint orbit through $K$.
Obviously, $\mathcal J_K$  is invariant under the adjoint action
of $Aut(\pi_{k,f})$ and the orbit space is a double quotient
space:
\begin{equation}\label{ orbit}
  Aut(\pi_{k,f}) \backslash  {\mathcal J_K} \cong  Aut(\pi_{k,f})
\backslash ( GL(3)/ G_K)
\end{equation}

 Now we give a scheme  to classify quadratic deformations of
Lie-Poisson structures  of Cases (7)-(10)  as follows.
\begin{itemize}
                    \item[(1)] Take a standard form of   Lie-Poisson
                    structure $\pi_{k,f}$
                    from the list in Theorem \ref{Th classify}(B) and take  a Jordan standard form $K$, for
                    which the
compatible 3-polynomial $F$ was  fixed  in \cite{LX}.
                     \item[(2)]  Choose a representative element
$\widetilde{K}$ in each  $Aut(\pi_{k,f})$ -orbit in ${\mathcal
J}_K$. This means that there exists a  $T \in GL(3)$ such that
   $\widetilde{K} = TKT^{-1}$. Then we get a compatible pair   $(\widetilde{K},
   \widetilde{F}) \in \mathcal{C}$ by Formula  (\ref{lx}) from the known compatible pair
   $(K, F)$.
                      \item[(3)] Check the    compatible pair   $(\widetilde{K},
   \widetilde{F})$  given above if  it satisfies  Equation (\ref{eq
                      1}). \item[(4)] In case that there are more than
one  cubic polynomials  satisfying  Equation (\ref{eq 1}), just
take one representative element.
\end{itemize}
By Theorem \ref{Th 2iso}, the Poisson structures $\pi_{k,f} +
\pi_{\widetilde{K}, \widetilde{F}}$ classify   quadratic
deformations of the Lie-Poisson structure $\pi_{k,f}$ after
checking for all Jordan standard forms. In following  sections we
classify quadratic deformations of Lie-Poisson structures.
Emphasizing again,  to save space we just write out some cases
with details. The others can be done by same way without any
difficulty except for some tedious computations.

\section{The  unimodular cases}

For  Lie-Poisson structures (1)-(6) listed in Theorem \ref{Th
classify},
      where their modular characters  vanish, $(K,F)\in \mathcal
      C_{0,f}$ if and only if $(K,F)\in \mathcal
      C$ and $K\in Der_0(\frak g)$ by Theorem \ref{Th deform}. The next theorem gives the forms of  $Der_0(\frak g)$, which are easy
      to be checked by Formula (\ref{eq
      tr0}).
\begin{thm}\label{Th der0}
       For Lie-Poisson structures (1)-(6), a pair (K,F) defines  a quadratic deformation  if and
only if $K$ has the following forms:
\begin{itemize}
            \item[]
$(1)  \quad  K \in {\frak sl(3)}, \, \, \, \quad  \quad (2) \quad
\, K \in{\frak o(3)}, \, \, \,  \quad  \quad (3)  \quad K \in
{\frak o(2, 1)},$
  \item[]$(4) ~~
                            \left(
                            \begin{array}
                            {ccc}
                                0           &  \alpha   & 0     \\
                                -\alpha     &   0       & 0   \\
                                   \beta    &  \gamma & 0

                            \end{array}
                            \right),
                             ~~~~~~~~~~
                             \quad (5)\, ~~ \left(
                            \begin{array}
                            {ccc}
                                0           &  \alpha   & 0      \\
                                \alpha      &   0       & 0     \\
                                   \beta   & \gamma   & 0
                            \end{array}
                            \right),
                             \quad  (6) ~~
                             \left(
                            \begin{array}
                            {ccc}
                                0           &      0      &     0 \\
                                  \alpha &  \delta  & \theta\\
                                   \beta   &  \gamma   & -\delta

                            \end{array}
                            \right).$
\end{itemize}
In all cases above, $F$ may be any cubic polynomial such that
$\hat{K}F =0$.
\end{thm}


 For Case (1) in Theorem \ref{Th classify}, the work  of the classification has been done in \cite{LX}. For other cases,
 One needs to choose firstly a representative element in $Der_0({\frak g})$ in each  adjoint orbit
                     of  $Aut({\frak g})$ given in Theorems \ref{Th
        Auto} and \ref{Th der0} respectively and then to find compatible homogeneous cubic polynomials. Here we only write out two examples  and leave
the others  to  interested readers.

        \noindent\textbf{$\bullet$ }   The compatible pair of ${\frak o(3)}$ is $(0,x^2+y^2+z^2)$ and  $Der_0{\frak o(3)}=Der{\frak o(3)}={\frak o(3)}$.
         Then by Theorem \ref{Th deform} all the
        quadratic deformations determined by those $(K,F)$  such that
        $K \in {\frak o(3)}$ and $\hat{K}F=0$. When $K=0$,
         $F$ may be any  cubic
        polynomial.  When $K\neq0$, by  Theorem \ref{Th Auto} and Theorem \ref{Th der0},
        any quadratic deformation is isomorphic to  one of the following forms:
        \begin{equation}
            K=  \left(
                \begin{array}{ccc}
                0           &  \alpha   &   0 \\
                -\alpha     &   0       &   0 \\
                0           &   0       &   0
                \end{array}
                \right),~~~~ F=a(x^2+y^2)z + bz^3, ~~~~ (\alpha >0, ~~ a\geq 0).
        \end{equation}
            In fact,  the automorphism group of ${\frak o(3)}$
            is $O(3)$ and any element of ${\frak o(3)}$ can be
            transformed to the above
            standard form via the adjoint action of $O(3)$.
            Moreover, the isotropy group of $K$ can change $z$ to
            $-z$ so that we can take $a\geq0$.

   \noindent\textbf{$\bullet$ } The compatible pair of ${\frak o(2,1)}$ is $(0,x^2+y^2-z^2)$.
Similar to the discussion above, one can check that when $K=0$,
         $F$ may be any  cubic
        polynomial and  when $K\neq0$,
        any
        quadratic deformation is isomorphic to  one of the following forms:
 \begin{itemize}
                            \item[\rm (1)]
             $ K_1=\left(
                \begin{array}{ccc}
                0           &  \alpha   &   0 \\
                -\alpha     &   0       &   0 \\
                0           &   0       &   0
                \end{array}
                \right), ~~~~ F_1=a(x^2+y^2)z + bz^3$
         \item[\rm (2)]
          $K_2=\pm\left(
            \begin{array}{ccc}
               0           &  0        &  1 \\
               0           &   0       &   1 \\
               1     &   1 &   0
            \end{array}
            \right), ~~~~~~~~F_2=a(x-y)^3 + b(x-y)(x^2+y^2-z^2)$
         \item[\rm (3)]
          $K_3=\left(
            \begin{array}{ccc}
            0           &   0       &   0        \\
            0           &   0       &   \alpha   \\
            0           &  \alpha   &   0
            \end{array}
            \right), ~~~~~~~~F_3=ax^3+ bx(y^2-z^2)$
 \end{itemize}
In fact,  ${\frak o(2,1)}$
        is isomorphic to ${\frak sl(2, \mathbb{R})}$  and any element of ${\frak o(2, 1)} $ can be
            corresponded to the one of above
            standard forms via the analysis the adjoint orbits of $SL(2,
            \mathbb{R})$.  The forms of $F_1$ and $F_3$ are obvious.
 $F_2$ can be checked after a coordinate transformation:
 $\xi=\frac{1}{2}(x+y), ~~\eta=\frac{1}{2}(x-y)$.

\section{Cases \, (7)-(10)}
In this section we follow the scheme shown  in Section 3. to
classify quadratic deformations of  other cases. First we consiter
the book algebra which is Case $(7)$ listed in Theorem \ref{Th
classify}. To save space we just write out details for three
Jordan standard forms only , given by (\ref{form K}), (\ref{form
K2}) and (\ref{form K3}) respectively.

Note that  the real projective plane $P^2$  is isomorphic  to a
$GL(3)$ homogenous space by dividing  an isotropy subgroup keeping
the subspace $\mathbb R e_3$ invariant. Such subgroup is just
$(G_7 \times \mathbb{R}^{\sharp}I)$, where $\mathbb
R^{\sharp}=\mathbb R-\{0\}$. That is,
 $$
P^2 \cong (G_7 \times \mathbb{R}^{\sharp}I) \backslash GL(3) \cong
G_7 \backslash GL(3) /\mathbb{R}^{\sharp}I,
$$
 Actually, such a correspondence is because any matrix of $G_7$ preserves $e_3$ then
the multiplication with it preserves the last column of any
matrix. Moreover, it is easy to see that $G_7G_K=G_KG_7$ and
$\mathbb{R}^{\sharp}I \subset G_K$ for any Jordan standard form
$K$. This means that, for Case $(7)$, the double orbit space
(\ref{ orbit}) has the following form:
\begin{equation}\label{p2} G_7 \backslash \mathcal
J_K \cong (G_7 \backslash GL(3))/G_K \cong (G_7 \times
\mathbb{R}^{\sharp}I) \backslash GL(3)/G_K
 \cong
P^2/G_K .
\end{equation}
For any  $[\alpha,\beta,\gamma]\in P^2$, where
$v=(\alpha,\beta,\gamma)$ satisfies $\parallel v\parallel=1$ and
$(v,e_3)\geq 0$, it corresponds an orthogonal matrix $T$. For $
v=e_3$, let $T = I$. For $( v,e_3)<1$, let
\begin{equation}\label{T}
w=(e_3-( v,e_3)v)/(\sqrt{1-( v,e_3)^2}),\quad T= (w, v\times w,
v)^T \in SO(3).
\end{equation}

   For a fixed Jordan standard form $K$ with its isotropy subgroup
 $G_{K}$ and a compatible pair $(K, F) \in \mathcal{C}$ given in \cite{LX}, let $p_i$ denote a  representative element of each orbit of $P^2$ under
the action of $G_{K}$ and $T_i$ is the corresponding  orthogonal
matrix given above. Then the pair $K_i=T_{i}KT_{i}^{-1}$ and
$F_i=F\circ T_{i}^{-1}$ is also a compatible pair on each  $G_{K}$
orbit by (\ref{lx}) (Here $det(T_{i})= 1$). First we take
\begin{equation}\label{form K}
K= diag(\lambda_1, \lambda_2, \lambda_3 ),\quad
\lambda_1+\lambda_2+\lambda_3=0, \quad\lambda_1\neq \lambda_2\neq
\lambda_3\neq 0.
\end{equation}
Its compatible cubic polynomials are in the form $F=axyz$ and
$G_K$ is the set of nonsingular diagonal matrices.
\begin{lem}\label{lemma3} For the  Jordan form $K$ given by (\ref{form K}), $P^2/G_K$ contains  seven orbits with following representative elements $p_i$
 and corresponding compatible pairs $(K_i,
 F_i)$:
\begin{itemize}
\item[(1)]$p_1=[0,0,1]$, \quad $K_1=diag(\lambda_1, \lambda_2,
\lambda_3 ),$\quad $F_1=axyz$.

\item[(2)]$p_2=[0,1,0]$,\quad$K_2=diag(\lambda_3, \lambda_1,
\lambda_2 ),$\quad $F_2=axyz$.

\item[(3)]$p_3=[1,0,0]$, \quad$K_3=diag(\lambda_3, \lambda_2,
\lambda_1 ),$\quad $F_3=-axyz$.
\item[(4)]$p_4=[\frac{\sqrt{2}}{2},\frac{\sqrt{2}}{2},0]$,~\quad$K_4= \frac{1}{2}\left(\begin{array} {ccc} 2\lambda_3  &  0 &0\\
0 & -\lambda_3& \lambda_1-\lambda_2
\\0 &\lambda_1-\lambda_2&-\lambda_3
\end{array}\right)$,~ \quad$F_4=\frac{a}{2}x(z^2-y^2)$.

\item[(5)]$p_5=[0,\frac{\sqrt{2}}{2},\frac{\sqrt{2}}{2}]$,\quad
$K_5=\left(\begin{array} {ccc} \frac{\lambda_2+\lambda_3}{2}  &  0
&\frac{\lambda_3-\lambda_2}{2}
\\
0 & \lambda_1 & 0
\\\frac{\lambda_3-\lambda_2}{2} &0&\frac{\lambda_2+\lambda_3}{2}
\end{array}\right)$,~ $F_5=\frac{a}{2}y(z^2-x^2)$.

\item[(6)]$p_6=[\frac{\sqrt{2}}{2},0,\frac{\sqrt{2}}{2}]$,\quad
$K_6=\left(\begin{array} {ccc} \frac{\lambda_1+\lambda_3}{2}  &  0
&\frac{\lambda_3-\lambda_1}{2}
\\
0 & \lambda_2 & 0
\\\frac{\lambda_3-\lambda_1}{2} &0&\frac{\lambda_1+\lambda_3}{2}
\end{array}\right)$,~ $F_6=\frac{a}{2}y(x^2-z^2)$.

\item[(7)]$p_7=[\frac{\sqrt{3}}{3},\frac{\sqrt{3}}{3},\frac{\sqrt{3}}{3}]$,\quad
$K_7=\left(\begin{array} {ccc}
\frac{\lambda_1+\lambda_2+4\lambda_3}{6}  &
\frac{\sqrt{3}(\lambda_2-\lambda_1)}{6}
&\frac{\sqrt{2}(2\lambda_3-\lambda_1-\lambda_2)}{6}
\\
\frac{\sqrt{3}(\lambda_2-\lambda_1)}{6} &
\frac{\lambda_1+\lambda_2}{2} &
\frac{\sqrt{6}(\lambda_1-\lambda_2)}{6}
\\\frac{\sqrt{2}(2\lambda_3-\lambda_1-\lambda_2)}{6} &\frac{\sqrt{6}(\lambda_1-\lambda_2)}{6}&\frac{\lambda_1+\lambda_2+\lambda_3}{3}
\end{array}\right)$,\\$F_7=\frac{\sqrt{3}a}{9}z^3+G(x,y,z)$, where $G(x,y,z)$ is cubic Polynomial without $z^3$.

\end{itemize}

\end{lem}
\pf  Note that for each $A= diag(m,n,p)\in G_{K}$,~ then
$[\alpha,\beta,\gamma]A=[m\alpha,n\beta,p\gamma]$. If
$\alpha=0,~\beta=0,~\gamma\neq0$, let $p=\frac{1}{\gamma}$,  we
obtain the first representative element $p_1$. The others are
same. $(K_i,
 F_i)$ can be get by (\ref{T}) easily.\qed
\begin{pro}\label{prop 1}
For  Lie-Poisson structure
$\pi=\frac{1}{2}\widehat{I}\wedge\widehat{k}$  ($k=(0,0,1)^T$) and
$K$ in form (\ref{form K}),  then any compatible pair
$(\widetilde{K},\widetilde{F}) \in \mathcal{C}_{k,0}$  such that
$\widetilde{K}\in \mathcal{J}_K$ is isomorphic to one of  Pairs
$(1),(2),(3)$ listed in Lemma \ref{lemma3} with $$ F_1 =
\frac{1}{6}(\lambda_2-\lambda_1)xyz, \quad F_2=
\frac{1}{6}(\lambda_1-\lambda_3)xyz, \quad
F_3=\frac{1}{6}(\lambda_3-\lambda_2)xyz.$$
\end{pro}
\pf  We only give the proof of Case (1) and Case (4) in Lemma
\ref{lemma3}
 and the proof of other cases is similar. First, for $(k,f)=((0,0,1)^T, 0)$ being the compatible pair of the book
 algebra,
Equation (\ref{eq 1}) in this case is:
\begin{equation}\label{eq 3}
 \frac{\partial F}{\partial
z}=-\frac{1}{6}(-a_{21}x^2+a_{12}y^2+(a_{11}-a_{22})xy+a_{13}yz-a_{23}xz),
\end{equation}
where $K= (a_{ij})$. Thus, for Case (1), one should has
$\frac{\partial F_1}{\partial
z}=\frac{1}{6}(\lambda_2-\lambda_1)xy$  by Formula (\ref{eq 3}).
On the other hand we have $F_1=axyz$, so if
$a=\frac{1}{6}(\lambda_2-\lambda_1)$, the equation is satisfied
and $(K_1,\frac{1}{6}(\lambda_2-\lambda_1)xyz)$ makes quadratic
deformation. As in Case (4), by Formula (\ref{eq 3}), we have
$$\frac{\partial F_4}{\partial
z}=-\frac{1}{6}[(\lambda_3-\frac{\lambda_1+\lambda_2}{2})xy-\frac{\lambda_1-\lambda_2}{2}xz].$$
We know that $F_4=\frac{a}{2}x(z^2-y^2)$ so that $\frac{\partial
F_4}{\partial z}=axz$. This implies that
$2\lambda_3=\lambda_1+\lambda_2$ $\Longrightarrow \lambda_3=0$
since $\lambda_1+\lambda_2+\lambda_3=0$, this is a contradiction
since we have $~\lambda_1\neq \lambda_2\neq \lambda_3\neq 0$, so
for any $a$,  Formula (\ref{eq 3}) couldn't be satisfied so that
there is no quadratic deformation in this case. \qed

Next we consider another Jordan form $K$ whose corresponding
compatible cubic polynomials(see \cite{LX}) and isotropy subgroup
are as follows:
\begin{equation}\label{form K2}
K=diag(\lambda, \lambda, -2\lambda),
\quad  F=mxyz+nx^2z+py^2z,  ~~~~~~\, ~~~~~~~~~\,  G_K=\left(\begin{array} {cc}A  &  0 \\
0 & a  \end{array}\right),
\end{equation}
where $\lambda\neq0$ and $A\in GL(2)$. It is easy to check that
$P^2/G_K$ has three orbits.
\begin{lem}\label{lemma4}With same notation in
Lemma \ref{lemma3} but $K$ being form (\ref{form K2}), then we
have
\begin{itemize}
\item[(1)]$p_1=[0,0,1]$,\quad $K_1=diag(\lambda, \lambda,
-2\lambda)$,\quad $F_1=mxyz+nx^2z+py^2z$.

\item[(2)]$p_2=[0,1,0]$,\quad $K_2=diag(-2\lambda, \lambda,
\lambda)$,\quad $F_2=mxyz+nxy^2+pxz^2$.

\item[(3)]$p_3=[0,\frac{\sqrt{2}}{2},\frac{\sqrt{2}}{2}]$,\quad$K_3=\left(\begin{array} {ccc} -\frac{1}{2}\lambda  &  0 &-\frac{3}{2}\lambda\\
0 &\lambda  & 0 \\-\frac{3}{2}\lambda &0& -\frac{1}{2}\lambda
\end{array}\right)$,\\[2mm] $F_3=\frac{m}{2}y(x^2-z^2)+\frac{\sqrt{2}}{2}ny^2(z+x)+\frac{\sqrt{2}}{4}p(z-x)^2(z+x)$.
\end{itemize}
\end{lem}
\begin{pro}\label{prop 2}
With same notation in Prop.\ref{prop 1} but $K$ being form
(\ref{form K2}),  then any compatible pair
$(\widetilde{K},\widetilde{F}) \in \mathcal{C}_{k,0}$  such that
$\widetilde{K}\in \mathcal{J}_K$ is isomorphic to one of  Pairs
$(1),(2)$ listed in Lemma \ref{lemma4} with $F_1 = 0$ and $ F_2=
\frac{1}{2}\lambda xyz+nxy^2$ respectively.
\end{pro}
 Finally we consider the following  Jordan form $K$ whose corresponding compatible cubic polynomials(see \cite{LX}) and
isotropy subgroup are as follows:
\begin{equation}\label{form K3}
 K=\left(\begin{array} {ccc} 0  &  1 &0\\
0 &0 & 1 \\0 &0& 0 \end{array}\right),
 \, \,F=pz^3+2qz^2x-qy^2z, ~~~~~~~\, ~~~~~~~~ \,  G_K=\left(\begin{array} {ccc} a  &  b &c\\
0 &a & b \\0 &0& a \end{array}\right).
\end{equation}
In this case it is easy to check that $P^2/G_K$ has three orbits
and two of them can make quadratic deformations.
\begin{lem}\label{lemma5}With same notation in
Lemma \ref{lemma3} but $K$ being form (\ref{form K3}), then we
have
\begin{itemize}
\item[(1)]$p_1=[0,0,1]$,\quad $K_1=\left(\begin{array} {ccc} 0  &  1 &0\\
0 &0 & 1 \\0 &0& 0 \end{array}\right)$,\quad
$F_1=pz^3+2qz^2x-qy^2z$.

\item[(2)]$p_2=[0,1,0]$,\quad  $K_2=\left(\begin{array} {ccc} 0  &  0 &0\\
0 &0 & 1 \\1 &0& 0 \end{array}\right)$,\quad
$F_2=px^3+2qx^2y-qz^2x$.

\item[(3)]$p_3=[1,0,0]$,\quad $K_3=\left(\begin{array} {ccc} 0  &  0 &0\\
-1 &0 & 0 \\0 &-1& 0 \end{array}\right)$,\quad
$F_3=px^3+2qx^2z-qy^2x$.
\end{itemize}
\end{lem}
\begin{pro}
With same notation in Prop.\ref{prop 1} but $K$ being form
(\ref{form K3}),  then any compatible pair
$(\widetilde{K},\widetilde{F}) \in \mathcal{C}_{k,0}$  such that
$\widetilde{K}\in \mathcal{J}_K$ is isomorphic to one  of Pairs
$(2),(3)$ in Lemma \ref{lemma5} with $F_2 =
px^3-\frac{1}{6}x^2y+\frac{1}{12}z^2x$ and $ F_3=
px^3-\frac{1}{6}x^2z+\frac{1}{12}y^2x $.
\end{pro}

For Lie-Poisson structures (8)-(10), by Theorem \ref{Th Auto}, we
know that $G_{l}\subset G_7 $  so that $G_{l}$-orbit~~ $\subset$
~~ $G_7$-orbit  $(l = 8, 9, 10)$ in $\mathcal{J}_K$. Therefore,
for these cases, the double orbit space (\ref{ orbit}) has the
following form by (\ref{p2}):
$$G_{l}\backslash \mathcal J_K     \cong  (G_{l}\backslash G_7)\cdot ( G_7 \backslash \mathcal J_K)
 \cong   (G_{l}\backslash G_7)  \cdot (P^2/G_K ). $$

It is easy to see that the quotient spaces   $G_{l}\backslash G_7$
$(l = 8, 9, 10)$  are 2-dim. manifolds and  their representative
matrices   can be given explicitly by fact that any invertible
matrix can be decomposed into product of an orthogonal matrix with
a lower triangular matrix (with a symmetric matrix for Case (9)).
For example,
$$G_{10}\backslash G_7\cong   S^1\times \mathbb{R}^{\sharp}  \cong \{ ~~T_{\alpha} Q_s = \left(\begin{array} {ccc} \cos\alpha  &  \sin\alpha &0\\
-\sin\alpha & \cos\alpha  & 0 \\0 &0& 1 \end{array}\right)\left(\begin{array} {ccc} 1 &  0 &0\\
0 & s & 0 \\0 &0& 1 \end{array}\right),
   ~ ~~  s\neq0\}.$$
Thus, for Case $(10)$, one has $G_{10}\backslash \mathcal J_K
 \cong   (S^1\times\mathbb{R}^{\sharp}) \cdot (P^2/G_K )$.
Now let $K$ be given in (\ref{form K}) and analyze  Case $(10)$.
From Lemma \ref{lemma3}, we know that $\mathcal J_K$ splits into
seven $G_7$-orbits, i.e.,
$$\mathcal J_K =  \bigcup_{i=1}^7 (G_7 \cdot K_i)   ~~ \Rightarrow ~~  G_{10} \backslash \mathcal J_K
 = \bigcup_{i=1}^7 G_{10}\backslash (G_7 \cdot K_i)
 \cong  \bigcup_{i=1}^7  (S^1\times\mathbb{R}^{\sharp}) \cdot K_i$$
Note that $ Q_s K_i = K_i Q_s$ and   $sF_i\circ Q_s^{-1} = F_i$
for    $i = 1, 2, 3, 5, 6$ ~~ in Lemma \ref{lemma3}.
 Consequently,
all $G_{10}$-orbits in $\mathcal J_K$, \,
$(S^1\times\mathbb{R}^{\sharp}) \cdot K_i$, \, can be
parameterized by $ i, \alpha$ and ($ i = 4, 7$) $ s$. Moreover, by
(\ref{lx}) and fact that $det(T_{\alpha}Q_s) = s$, we can get
their representative compatible pairs:
$$(  T_{\alpha}K_iT_{\alpha}^{-1} ,~~\, F_i\circ T_{\alpha}^{-1}), ~~ \, \,  i = 1, 2, 3, 5, 6, ~~ \,  ~~ \, \alpha \in [0, 2\pi)$$
$$(  (T_{\alpha}Q_s)K_i(T_{\alpha}Q_s)^{-1} ,~~\, sF_i\circ(T_{\alpha}Q_s)^{-1}), ~~ \, \,  i = 4, 7, ~~ \, s\neq 0,  ~~ \, \alpha \in [0, 2\pi).$$
The next step is to figure out all  representative compatible
pairs   satisfying  Equation (\ref{eq 1}) from those given above
to classify quadratic deformations of Lie-Poisson structure Case
 (10) on  $\mathcal J_K$. It is easy to check that  there is no  quadratic deformation
 for $ i = 1, 2, 3$ listed in Lemma \ref{lemma3}.

As the  last example, we  take $K$  as given in (\ref{form K2}).
In this case we know that $\mathcal J_K$ splits into three
$G_7$-orbits
  with their representative compatible
pairs $( K_i , F_i), ~~ i = 1, 2, 3$ listed in  Lemma
\ref{lemma4}.  Here we only write out the conclusion for $K_1 =
K$.
\begin{pro} For Lie-Poisson structure Case (10) and  $K$ being
(\ref{form K2}),
 then any  quadratic deformation
$(\widetilde{K},\widetilde{F})$  such that $\widetilde{K}\in
G_7\cdot K$ is isomorphic to $( K, ~~-2\lambda x^2z).$
\end{pro}
\pf It is obviously that $G_7\cdot K=G_{10}\cdot K$ in this case
since $ (T_{\alpha}Q_s)K(T_{\alpha}Q_s)^{-1}= K $ so that  we can
take $( K_1, F_1)$ given in Lemma \ref{lemma4} for $i = 1$ as the
representative compatible pair on this orbit, which satisfies
Equation (\ref{eq 1}) if and only if $\frac{\partial F_1}{\partial
z}=-2\lambda x^2$. Thus when $m=p=0, ~n=-2\lambda$, the equation
is satisfied so that this pair makes quadratic deformation.\qed

\end{document}